\documentclass[12pt]{article}

\usepackage{fancyhdr, amsmath, amssymb, amsfonts, url, subfigure, dsfont, mathrsfs, graphicx, epsfig, subfigure, amsthm, tikz}

\usepackage{todonotes}

    \usepackage{epsfig}
    \usepackage{graphicx}
    \usepackage{color}
    \usepackage{wasysym}
 
\newlength{\guillotine}
\newtheorem{thm}{Theorem}[section]
\newtheorem{cor}[thm]{Corollary}
\newtheorem{lemma}[thm]{Lemma}

\newtheorem{definition}[thm]{Definition}

\theoremstyle{remark}
\newtheorem{rem}[thm]{Remark}
\begin{document}

\title{Accurate  bounds on Lyapunov exponents \\
for expanding maps of the interval}

\author{M. Pollicott and P. Vytnova\thanks{ The first author is partly
supported by ERC-Advanced Grant 833802-Resonances and EPSRC grant
EP/T001674/1 the second author is partly supported by
EPSRC grant EP/T001674/1. }}

\date{ }

\maketitle

\abstract{In this short note we describe a simple but remarkably effective
method for rigorously estimating Lyapunov exponents for expanding maps of the
interval.  We illustrate  the applicability of this method with some standard examples.} 


\section{Introduction}
\label{s:intro}
Lyapunov exponents  give a well known characterization of the instability in a dynamical system by quantifying how nearby orbits separate.
In particular, a non-zero Lyapunov exponent with respect to an invariant ergodic measure
implies that typical nearby  orbits separate exponentially quickly. It is therefore useful
to have a rigorous and effective estimate of these values, in particular, 
in the setting of one dimensional expanding maps for an absolutely continuous
invariant probability measure. This problem has attracted the attention of many authors 
who have employed a variety of different methods (see~\cite{bose},~\cite{keane},~\cite{wormell}).


In this paper we will consider the nice class~$NC$ of expanding piecewise analytic mixing Markov maps of the
interval. We recall the definition. 
\begin{definition}
    \label{def:markovMap}
Let $I = [a,b]$ be a closed interval. We say that a map~$f$ belongs to the
class~$NC(I)$ if  there exists a partition $a=x_1 < x_2 < \cdots < x_{n+1} = b$ such that:
 \begin{enumerate}
     \item The restrictions $f|_{[x_j, x_{j+1}]}$ are analytic maps for $j = 1, \ldots n$;
     \item There exists $d>1$ such that for all~$1\le k \le n$ and for all $ x \in (x_k,
         x_{k+1})$ we have $|f'(x)| \geq d$. 
%
     \item The Markov property holds: if $f((x_k,x_{k+1})) \cap (x_j,x_{j+1}) \neq \varnothing$;
         then $f( (x_k,x_{k+1})) \supset (x_j,x_{j+1})$.
         \item 
         The map $f$ is topologically mixing\footnote{This is equivalent to the
         map being locally eventually onto, i.e. to saying that there exists~$N \geq 1$ such that for each $1\leq j \leq n$ we have
     $\overline{f^N( (x_j,x_{j+1}) )} = I$.}
 (i.e., for any non-empty open sets $U,V \subset I$ there exists $n_0\geq 1$
 such the for all $n \geq n_0$ we have $U\cap f^{-n}V \neq \varnothing$).
 \end{enumerate}
 \end{definition}

For every map~$f \in NC(I)$ there exists a unique absolutely continuous~$f$-invariant 
probability measure $d\mu = \rho(x) dx$ on~$I$ \cite{collet}.  In particular, the measure~$\mu$ is ergodic.
Furthermore, every map of the class $NC(I)$ is invertible on each of the intervals
$(x_j,x_{j+1})$, in particular, there exist analytic maps $f_{jk}\colon (x_j, x_{j+1}) \to (x_k, x_{k+1})$,  
such that $f (f_{jk}(x)) = x$
whenever  $f((x_k,x_{k+1})) \cap (x_j,x_{j+1}) \neq \varnothing$.
The maps $f_{jk}$ are called \emph{inverse branches} of~$f$.

A standard  approach  to constructing  the  
measure~$\mu$  is to use transfer operators. 
Let us denote by~$\mathcal B$ the space of analytic functions on the disjoint union $\coprod_{k=1}^n[x_k,x_{k+1}]$.
We can introduce a one-parameter family of linear operators 
$\mathcal L_t: \mathcal B \to \mathcal B$ 
($t \in \mathbb R$)
called {\it transfer operators} defined in terms of inverse
branches of~$f$:  
\begin{equation}
    \label{eq:trop}
[\mathcal L_t h](x) = \sum_{k: f((x_k,x_{k+1})) \cap (x_j,x_{j+1}) \neq
\varnothing} |f_{jk}^\prime(x)|^t h(f_{jk}(x))\chi_{[x_j, x_{j+1}]}(x),
\end{equation}
where $\chi_{[x_j, x_{j+1}]}$ is the indicator function of the interval $[x_j,
x_{j+1}]$. 
 In the special case that 
 $f$ is full branched, i.e., 
 $f((x_k,x_{k+1})) = I$ for $k=1, \cdots, n$,  we can denote the inverses
$f_k: I \to (x_k, x_{k+1})$, i.e.,   $f(f_k(x))=x$ for all $a \leq x \leq b$.  
The transfer operators 
$\mathcal L_t: C^\omega(I)  \to C^\omega(I) $ ($t \in \mathbb R$)  then  take the form
\begin{equation}
    \label{eq:trop2}
[\mathcal L_t h](x) = \sum_{j=1}^n |f_j^\prime(x)|^t h(f_j(x)), \qquad x \in I.
\end{equation}
 Most of our examples will be of this type.
 
It is well known~\cite{collet} that  the positive density~$\rho \in C^\omega(I)$ 
of the measure~$\mu$ is characterized as a fixed point for the operator
$\mathcal L_1$, corresponding to the parameter choice $t=1$.  
Nevertheless, including this operator into a one parameter family will serve us well later. 

We can now define the Lyapunov exponent of the system $(I,f,\mu)$ which quantifies the sensitivity
of typical orbits on initial conditions.  

\begin{definition} 
    \label{def:lyapexp}
We define the {\it Lyapunov exponent} for the map~$f$ and its stationary measure~$\mu$ by 
$$
\lambda(f, \mu) = \int_I \log |f'(x)| d\mu(x).
$$
\end{definition}

\begin{rem}
This value coincides with the metric entropy~$h(\mu)$ of the measure~$\mu$ by
the Rokhlin's formula \cite{yuri}.
\end{rem}

Since the measure~$\mu$ is ergodic,  
applying the Birkhoff ergodic theorem one can see that for $\mu$-almost all
$x\in I$ we get
$$
\lim_{n \to +\infty} \frac{1}{n} \log |(f^n)'(x)| = \lambda(f, \mu) .
$$
There are various methods used to estimate the Lyapunov exponents.
Probably the most famous are  Ulam's method and finite section
methods~\cite{keane}. Another approach is based on periodic points method~\cite{jenkinson}. 
Recent work by Wormell~\cite{wormell} is based on the Galerkin spectral method originally
developed for PDEs. 
In this note we present an alternative approach, which starts with the
\emph{spectral Chebyshev
collocation method}, also initially developed for PDEs~\cite[\S3]{canuto}. In
dynamical systems it has been
used succsefully by Babenko and Yuriev in their solution of the Gauss
problem~\cite{babenko79} and by Babenko in his computation of the fixed point of
the renormalisation operator for the period-doubling map~\cite{babenko83}.  In
our approach, we combine the Chebyshev collocation method with a small amount
of thermodynamic formalism (involving the pressure function) and a classical
min-max method. The  main advantage of this combination of ideas  is that it provides an efficient and 
effective way to estimate Lyapunov exponents and gives rigorous estimates with validated error bounds.

The main results  we present in this note are  the following.
The first theorem gives a method for obtaining rigorous bounds on the Lyapunov exponent.
\begin{thm}
    \label{thm:main}
Let $f\colon I\to I$ be an expanding piecewise analytic mixing Markov map of the interval with absolutely continuous probability measure $\mu$.
Assume that for some $\varepsilon > 0$ there exists a pair of positive
functions\footnote{The reason for this choice of notation is that in practice
the functions $p$ and $q$ are polynomials.}~$p, q \colon I \to \mathbb R^{+}$ and a pair of numbers
$0 < \alpha < \beta$ such that 
\begin{equation}
    \label{eq:lyaplims}
\sup_I\frac{ \mathcal L_{1+ \varepsilon} p}{ p } \leq  e^{-\alpha}
\hbox{ and } 
 \sup_I\frac{
\mathcal L_{1- \varepsilon} q}{ q } \leq  e^{\beta}.
\end{equation}
Then the following double inequality holds:
$$
 \frac{\alpha }{\varepsilon} \leq \lambda(f, \mu) \leq \frac{\beta
}\varepsilon. 
$$
\end{thm}

\begin{rem}
    The idea behind  Theorem~\ref{thm:main} is that for the class of maps we
    consider for any positive function~$p$ the supremum of the ratio $\frac{\mathcal
    L_tp}{p}$ gives an upper bound on the leading eigenvalue of~$\mathcal L_t$. 

    Note that if the function~$p$ is close to the leading eigenfunction
    of the operator~$\mathcal L_{t}$
    then the ratio $\frac{\mathcal L_{t}  p}{p}$ is close to a
    constant function. This observation allows us to estimate the ratios
    rigorously in practice.  
\end{rem}

\begin{rem}
If we do not assume that~$f$ is Markov then the statement of the
Theorem remains true, however, in this setting the construction of the
functions~$p$ and~$q$ is more challenging since the eigenfunctions of $\mathcal
L_t$ might be non-analytic (but of bounded variation). As we will see later, 
in practical applications, the interval $(\alpha,\beta) \ni \lambda(f,\mu)$
depends on the quality of approximation of the leading eigenfunction of~$\mathcal L_t$
by polynomials~$p$ and~$q$. 
\end{rem}

The next theorem guarantees that the previous theorem can be used to get bounds
on the Lyapunov  exponent which are arbitrary accurate. 
Note that Theorem~\ref{thm:main} also holds under the weaker assumption that 
$f\colon I\to I$ is an expanding piecewise~$C^2$ mixing Markov map of the
interval, however, in this case it is much harder to compute the functions~$p$
and~$q$ which will give us good estimates on the Lyapunov exponent.  
In addition, it is convenient to assume analyticity in order to apply the following theorem.
\begin{thm}
    \label{thm:justification}
Let $f\colon I\to I$ be an expanding piecewise analytic mixing Markov map of the interval with absolutely continuous probability measure $\mu$.
Then for any $\delta > 0$ we can choose $\varepsilon > 0$, $0 < \alpha  < \beta$ and strictly positive 
polynomials $p, q : I \to \mathbb R$ satisfying~\eqref{eq:lyaplims}
with 
\begin{equation}
\left|\frac{\beta}{\varepsilon} - \frac{\alpha}{\varepsilon}\right| < \delta.
\end{equation}
\end{thm}

\section{Examples}
\label{s:examples}
In this section we will demonstrate how Theorem~\ref{thm:main} can be used in
practice. To this end we consider four examples, and compare the 
estimates we obtain for the Lyapunov exponents with previously  known results.
 
Theorem~\ref{thm:main} allows us to obtain rigorous bounds using  
the built-in {\tt MaxValue} routine in Mathematica, and the implementation is 
relatively straightforward. However, some care is required in choosing
parameters during the construction of the functions~$p$ and~$q$. 
In Section~\ref{ss:justify} we give more details on the practicalities of the implementation. 

\subsection{Classical example: the Lanford map}
    \label{ex:lanford}
\begin{figure}[h!]
 \centerline{ \includegraphics[height=5.0cm,angle=0 ]{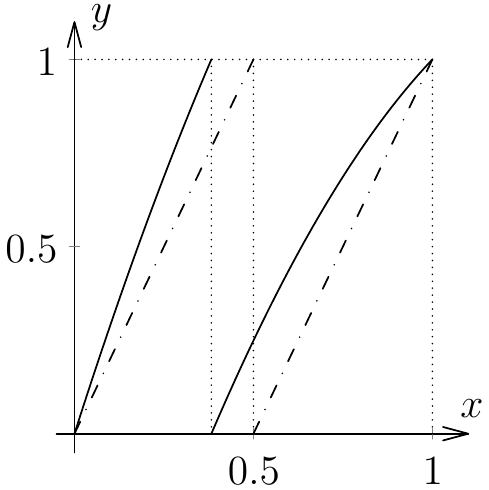}  }
 \caption{The solid curve is a  plot of the Lanford map.  For comparison, the
 dashed line is the plot of the doubling map.}
 \label{fig:lanford}
 \end{figure}

We  will first  illustrate our approach with  the standard example of the Lanford map~\cite{lanford}.
The original Lanford map $f_L:[0,1] \to [0,1]$ is  defined by
\begin{equation}
    \label{eq:lanford}
f_L(x) = 2x + \frac{1}{2}x(1-x) \quad \mod 1
\end{equation}  
and the graph of $f_L$ is shown in Figure~\ref{fig:lanford}.
  Observe  that the map is uniformly expanding   with 
  $T'(x) \geq  T'(1) = \frac{3}{2}$ for $0 \leq x \leq 1$.
  
 The inverse branches of $f_L$ are contractions given by
$$
f_1(x) = \frac{5 - \sqrt{25 - 8x}}{2}
\quad \mbox{ and } \quad f_2(x) = \frac{5 - \sqrt{17 - 8x}}{2}.
$$
The transfer operator $\mathcal L_t$ therefore takes the form
$$
[\mathcal L_th](x) = 
\left(\frac{2}{ \sqrt{25 - 8x}}\right)^t h\left(\frac{5 - \sqrt{25 - 8x}}{2} \right)  +  
\left(\frac{2}{\sqrt{17 - 8x}}\right)^t h\left(\frac{5 - \sqrt{17 - 8x}}{2} \right).
$$
Due to simplicity of the formulae involved we shall attempt to obtain estimates
on the Lyapunov exponents of particularly high accuracy, to demonstrate the
power of our method. Namely, we shall choose $\varepsilon = 10^{-180}$. Then we
fix $N=400$ and compute the nodes of the Chebyshev polynomial $T_{400}$ to $512$-digits
precision. Subsequently, we want to construct the functions~$p$ and~$q$ as polynomials of degree $399$
using the spectral Chebyshev collocation method. 
We validate that the polynomials $p$ and $q$ are positive using the method
described in Section~\ref{ss:justify}. At this point we apply the built-in {\tt
MaxValue} Mathematica routine to calculate
$$
\alpha: = -\log \left( \mbox{\tt MaxValue} \frac{\mathcal L_{1+\varepsilon}
p}{p}\right) \mbox{ and } 
\beta: = \log \left( \mbox{\tt MaxValue} \frac{\mathcal L_{1-\varepsilon} q}{q}
\right).
$$
We obtain the following values: 
\begin{align*}
\alpha = 6.&5766178000\,6597677541\,5824138238\,3206574324\,1069580012\,2019539528\, \\ 
               &0269163266\,6111554023\,7595564597\,5291517482\,9642156331\,7980263014\, \\ 
               &8859489891 \times 10^{-181} 
               \mbox{ and }
\end{align*}
\begin{align*}
\beta = 6.&5766178000\,6597677541\,5824138238\,3206574324\,1069580012\,2019539528\, \\
              &0269163266\,6111554023\,7595564597\,5291517482\,9642156331\,7980263014\, \\ 
              &8859489094 \times 10^{-181}; 
\end{align*}
which are each presented to $130$ significant figures.
In particular, with these choices 
Theorem~\ref{thm:main} yeilds that $\varepsilon^{-1} \alpha < \lambda(f_L,\mu) <
\varepsilon^{-1}\beta$ therefore 
$$
\begin{aligned}
\lambda(f_L, \mu)  =
0.&65766178000\,6597677541\,5824138238\,3206574324\,1069580012\,2019539528\, \\ &
  0269163266\,6111554023\,7595564597\,5291517482\,9642156331\,7980263014\, \\ & 
  88594891 \pm 10^{-128}\cr
\end{aligned}
$$
This value has previously been computed by Wormell~\cite{wormell} and her result
agrees with the above.
In the present approach, the simplicity of the functions
$p$ and $q$ is the source of the efficiency of the approach. In particular, this estimate was obtained
in approximately 2 hours on a \emph{personal} Macbook pro laptop with 2.8 GHz Quad-Core Intel Core i7 and
16 GB 2133 MHz LPDDR3 using Mathematica.


\begin{rem}
In addition to using the internal {\tt MaxValue} function,  whose  code is not
available to the public, we can apply a simple Monte-Carlo type method to numerically  verify the
value we obtained. More precisely, we generate $N_{mc}=1000$ pseudo-random points
$x_j$, $j=1, \ldots, 1000$ in the interval $[0,1]$ and evaluate both ratios at
these points to get the values 
$$
y_j^{+}: = \frac{[\mathcal L_{1+\varepsilon} p](x_j)}{p(x_j)} \qquad \mbox{ and
} \qquad 
y_j^{-}: = \frac{[\mathcal L_{1-\varepsilon} q](x_j)}{q(x_j)}, \quad j =
1,\ldots,100. 
$$
Then we compute $a_1:=-\log\max_j y_j^{+}$ and
$b_1:=\log \max_j y_j^{-}$. Repeating this procedure 
a total of  $t_{mc}=100$ times, we obtain the values that are within  a  distance of
$10^{-345}$ from~$\alpha$ and~$\beta$, respectively.
In particular,   we see that our estimate agrees with the estimate given by the function {\tt MaxValue}.   
\end{rem}

\subsection{Lanford family of maps}
\label{ss:lanfamily}
We can extend the first example by including it in a larger family of maps.
More precisely, 
we can include the Lanford map~\eqref{eq:lanford} into a family of expanding
maps $f_c:[0,1] \mapsto [0,1]$ defined by 
$$
f_c(x) \colon = 2x + c x(1-x) \quad  \mod 1, \qquad 0 < c <1. 
$$  
Observe that for the chosen parameter values 
$f_c'(x) =   2 + c (1-2x)  \geq 2 - c = f_c'(1)$ for $0 \leq x \leq 1$ and so the map $f_c$ is expanding.
Then the inverse branches  $f_1, f_2 \colon I \to I$ are contractions defined by
$$
f_1(x) = \frac{2+ c - \sqrt{(2+c)^2 - 4 c  x}}{2 c }
\quad \mbox{ and } \quad f_2(x) = \frac{2+ c - \sqrt{(2+ c)^2 - 4 c  (x+1)}}{2c }.
$$
Following the formula~\eqref{eq:trop} we obtain the associated  transfer operator
$\mathcal L_t$:
$$
\begin{aligned}
\left[\mathcal L_t h\right]  (x) =& 
\left(\frac{1}{\sqrt{(2+ c )^2 - 4  c  x}}\right)^t \cdot
 h\left(\frac{2 + c   - \sqrt{(2+ c )^2- 4  c  x}}{2 c } \right)  \cr
 &\qquad +  
\left(\frac{1}{\sqrt{(2+ c )^2 - 4  c  (x+1)}}\right)^t \cdot
 h\left(\frac{2 + c   - \sqrt{(2+ c )^2- 4  c  (x+1)}}{2 c } \right).
 \end{aligned}
$$
We next want to compute the Lyapunov exponent
$\lambda(c):=\lambda(f_c,\mu_c)$ for forty equally spaced values 
$c = c_j = 0.001+\frac{j-1}{40}$, with $j = 1,\dots,40$ with an error of~$10^{-3}$ to
sketch a graph of~$\lambda$ as a function of~$c$. For this purpose we choose
$\varepsilon=10^{-3}$ and $m=60$ and compute the nodes of the Chebyshev
polynomial $T_{60}$ with accuracy of $256$ digits. We then apply 
Theorem~\ref{thm:main} and obtain lower and upper
bounds for the Lyapunov exponent. The precision of the {\tt MaxValue} routine in
the computation was set to~$128$ digits.  

\begin{figure}[h!]
 \centerline{ \includegraphics{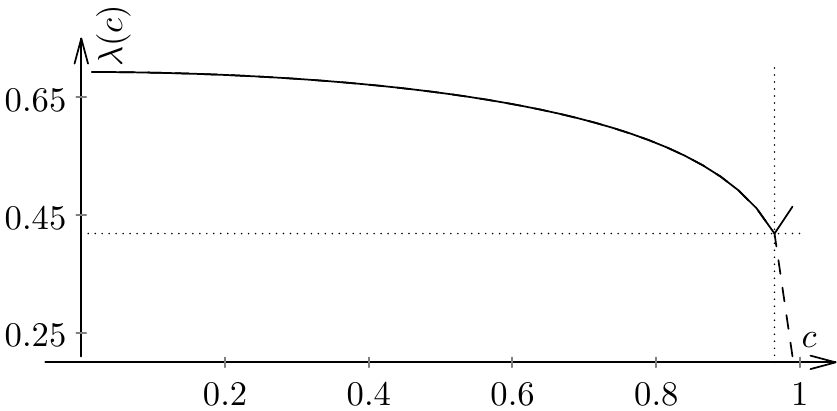}  }
 \caption{A plot of the Lyapunov exponent for the Lanford family 
 based on the calculation for the $40$ parameter values in
 the interval $[0.001,0.99]$. The dependence of the Lyapunov exponent on~$c$ is analytic. It 
  appears that the derivative $\frac{d \lambda(f_c,\mu_c)}{d c} \to -\infty$
  as~$c \to 1$. At the other end we have $\lambda(f_c,\mu_c) \to \log 2 =
 0.693\ldots$  as $c \to 0$ which is expected, since $c=0$ corresponds to the doubling map. }
 \label{fig:lanfordCurve}
 \end{figure}

Based on this calculation, we sketch the functions
$\frac{\alpha(c)}{\varepsilon}$ (dashed curve) and $\frac{\beta(c)}{\varepsilon}$
(solid curve) in Figure~\ref{fig:lanfordCurve}. We see that for $0.01 < c <
0.96$ the two curves are indistinguishable. However in the interval $ 0.96 < c
<0.99$ they appear to be different. This reflects the fact that 
$f_c^\prime(1)\to 1$ as~$c\to 1$, i.e. the map $f_c$ has weak hyperbolicity
for~$c$ close to~$1$. Uniform hyperbolicity is essential for
Theorems~\ref{thm:main} and~\ref{thm:justification} to be applicable. 

In addition, we may  also calculate the Lyapunov exponent for a selected parameter value
$c = \frac{1}{4}$, for example,  with  high  accuracy. To this end, we choose
$\varepsilon=10^{-180}$ and compute~$300$ zeros of the Chebyshev
polynomial~$T_{300}$ with accuracy of~$400$ digits. Then we apply the spectral
collocation method to construct polynomials~$p$ and~$q$ of degree~$299$.  
As before, we verify this this functions are positive, and apply {\tt MaxValue}
with working precision~$400$. 
$$
\alpha: = -\log \left( \mbox{\tt MaxValue} \frac{\mathcal L_{1+\varepsilon}
p}{p}\right) \mbox{ and } 
\beta: = \log \left( \mbox{\tt MaxValue} \frac{\mathcal L_{1-\varepsilon} q}{q}
\right).
$$
We obtain the following values (for which we give $166$ digits): 
\begin{align*}
\alpha=0.&6851020685\,7610906837\,8941120635\,3368474791\,2954208389\,7263352003\,\\
         &7686275679\,0996831645\,2222918013\,3822749913\,1527755618\,1523970004\,\\
         &1829353798\,5819153203\,8804954205\,2390123411\,591687 \times 10^{-180}; \mbox{ and } \\
\beta=0.&6851020685\,7610906837\,8941120635\,3368474791\,2954208389\,7263352003\,\\
        &7686275679\,0996831645\,2222918013\,3822749913\,1527755618\,1523970004\,\\
        &1829353798\,5819153203\,8804954205\,2390123411\,591699 \times 10^{-180}.
\end{align*}
This gives the value of the Lyapunov exponent with accuracy of $164$ decimal
places: 
\begin{align*}
\lambda\left(f_{\frac14},\mu_{\frac14}\right) &= 
      0.6851020685\,7610906837\,8941120635\,3368474791\,2954208389\,7263352003\,\\
        &7686275679\,0996831645\,2222918013\,3822749913\,1527755618\,1523970004\,\\
        &1829353798\,5819153203\,8804954205\,2390123411\,59169\pm10^{-165}.
\end{align*}
Using the Monte-Carlo method with $N_{mc}=1000$ pseudo-random points in the interval
$[0,1]$ and $t_{mc}=100$ samples, we can numerically  check  the output of the routine {\tt
MaxValue}. Namely, taking the maximum of the ratios $\frac{\mathcal
L_{1+\varepsilon} p}p$ and $\frac{\mathcal L_{1-\varepsilon}q}q$ computed at $1000$ different points a
hundred times, we obtain the values which lie within the distance of 
$2.0\times 10^{-345}$ from $\alpha$ or $\beta$, respectively. 

\subsection{A family of full branch piecewise M\"obius maps}
We next consider a family of  examples studied by Slipantschuk, Bandtlow and Just in~\cite{BJS} 
in connection with their study of relation between Lyapunov exponents and mixing
rates. 

Following~\cite{BJS}, for $- \frac{1}{4} \leq c \leq \frac{1}{2}$  we have a  map $f_c \colon [-1,1] \to [-1,1]$ defined by 
$$
f_c(x) = \frac{1-2(c+1)|x|}{1 + 2 c |x|}.
$$
When $c=0$ this reduces to a piecewise linear ``tent map''.  
In the special case $c=0.11$, of particular  importance  to the authors of~\cite{BJS},
they assert that the Lyapunov exponent is $\lambda(f_{0.11},\mu_{0.11}) = 0.685 \ldots$,
although  the paper does not provide any details as to how this value was computed. 

\begin{figure}
    \begin{center}
    \begin{tabular}{ccc}
 \includegraphics{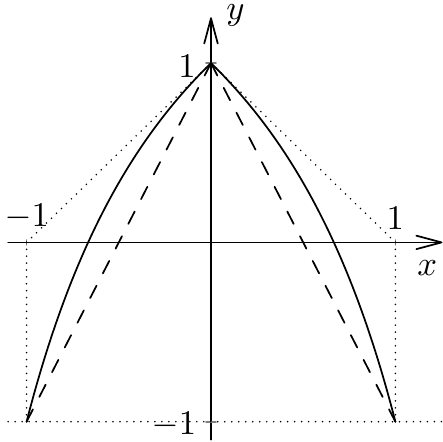} & & \includegraphics{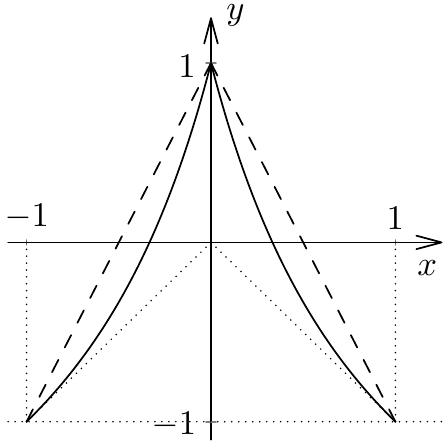} \\
 (a) & \phantom{Plot of the tent map}  & (b) \\
 \end{tabular}
 \end{center}
 \caption{Two plots of the bent tent map for parameter values close to the ends
 of the parameter interval: $c = 0.2495$ (a) and $c = 0.4995$ (b). The dashed 
 lines are the tent map corresponding to~$c=0$.}
 \end{figure}

\begin{figure}
    \centerline{\includegraphics{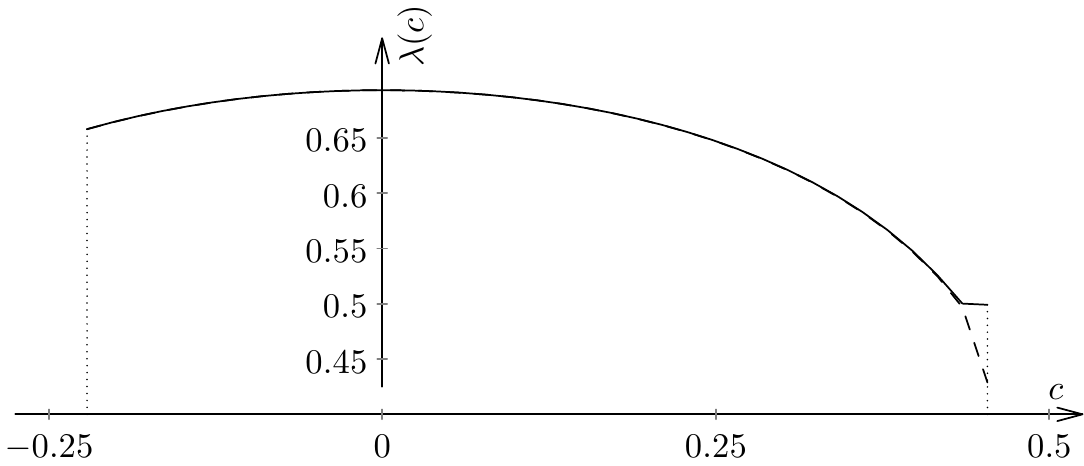}}
    \caption{Lower (dashed curve) and upper (solid curve) bounds on the Lyapunov exponent 
    for the family of bent  tent maps. We see that for $c \in (-0.24, 0.42)$ they are
    almost indisitnguishable. This is the range of parameter values where our
    method is particularly effective. }
    \label{fig:tentCurve}
\end{figure}

The inverse branches $f_1, f_2: [-1,1] \to [-1,1]$ take the form 
$$
f_1(x) = \frac{1-x}{2cx + 2(c+1)} \qquad \mbox{ and } \qquad f_2(x) = - \frac{1-x}{2cx + 2(c+1)}.
$$
In particular, $|f_1'(x)| = |f_2'(x)| = \frac{1+2c}{2(1 + c + x )^2}$.
The associated transfer operator is given by 
\begin{multline*}
\left[\mathcal L_t h \right](x) = \left|\frac{2 (2c+1)}{(2c x + 2(c+1)^2)}\right|^t h\left( \frac{1-x}{2cx + 2(c+1)} \right)  \\ + 
\left|\frac{2 (2c+1)}{(2c x + 2(c+1)^2)}\right|^t h\left(- \frac{1-x}{2cx + 2(c+1)} \right).
\end{multline*}
We shall   recover and  improve the estimate of~\cite{BJS}.
For this purpose, we choose $m=400$ and compute Chebyshev nodes with accuracy of
$600$ digits. Then we choose $\varepsilon = 10^{-175}$ and apply Chebyshev
collocation method to obtain two polynomials~$p$ and~$q$ of degree~$399$. We
then verify that they are positive and evaluate 
$$
\alpha: = -\log \left( \mbox{\tt MaxValue} \frac{\mathcal L_{1+\varepsilon}
p}{p}\right) \mbox{ and } 
\beta: = \log \left( \mbox{\tt MaxValue} \frac{\mathcal L_{1-\varepsilon} q}{q}
\right).
$$
with working precision set to~$400$. For each of the values we give~$180$
digits. 
\begin{align*}
\alpha = 0.&6849333272\, 2256432968\, 5622546648\, 2230532357\, 7867689297\, 3987148578\,   
         \\&  8085505250\, 5345328689\, 5040861069\, 9964717724\, 0662692746\, 4804164759\,
         \\&  1723161867\, 2782003116\, 7550103160\, 3289137884\, 1128687391\,
         8360864512 \times 10^{-175};  \\ 
\beta =  0.&6849333272\, 2256432968\, 5622546648\, 2230532357\, 7867689297\, 3987148578\, 
         \\&   8085505250\, 5345328689\, 5040861069\, 9964717724\, 0662692746\, 4804164759\,
         \\&   1723161867\, 2782003116\, 7550103160\, 3289137884\, 1128687391\, 8360866430 \times 10^{-175}.
\end{align*}
This yields the following estimate on the Lyapunov exponent accurate to $176$
decimal places given below:
$$
\begin{aligned}
    \lambda(f_{0.11}&,\mu_{0.11})  = 0.6849333272\, 2256432968\, 5622546648\, 2230532357\, 7867689297\, 3987148578\, \\ 
           & 8085505250\, 5345328689\, 5040861069\, 9964717724\, 0662692746\, 4804164759\, \\
           & 1723161867\, 2782003116\, 7550103160\, 3289137884\, 1128687391\,
           8360865 \pm  10^{-177}. 
\end{aligned}
$$
In addition, similarly to the case of the Lanford map, we can plot the Lyapunov exponent as
a function of the parameter~$c$. A sketch of the graph $\lambda(c)$ is shown in
Figure~\ref{fig:tentCurve}. It is based on the computation for $40$ equidistant
points in the parameter interval $(-0.24,0.45)$. The following setup
has been used for the calculation: $\varepsilon = 0.001$, $N = 128$ Chebyshev nodes computed with
accuracy of~$512$ digits. For the parameter values $c \in (-0.25,-0.24)$ and $c
\in (0.45, 0.5)$ the computation turns to be unstable and the resulting values
of $\alpha$ and $\beta$ disagree by as much as~$1.3$ for $c = -0.22$. This is
again due to the fact that $|f_c^\prime(0)| \to 1$ as $c \to -0.25$ and
$|f_c^\prime(\pm 1)|\to 1$ as $c \to 0.5$, i.e. diminuishing hyperbolicity of
the system.

\subsection{Bent baker's map}
Finally, we consider an example studied by Froyland in~\cite{froyland}.
Namely, we can consider the map~$f\colon [0,1] \to  [0,1]$ defined by
$$
f(x) = \frac{4\sqrt{6}}{3}x^3 - 2 \sqrt{6} x^2 + \left( 2 + \frac{2\sqrt{6}}{3}\right)x \quad
 \mbox{mod $1$.}
$$

\begin{figure}
 \centerline{ \includegraphics{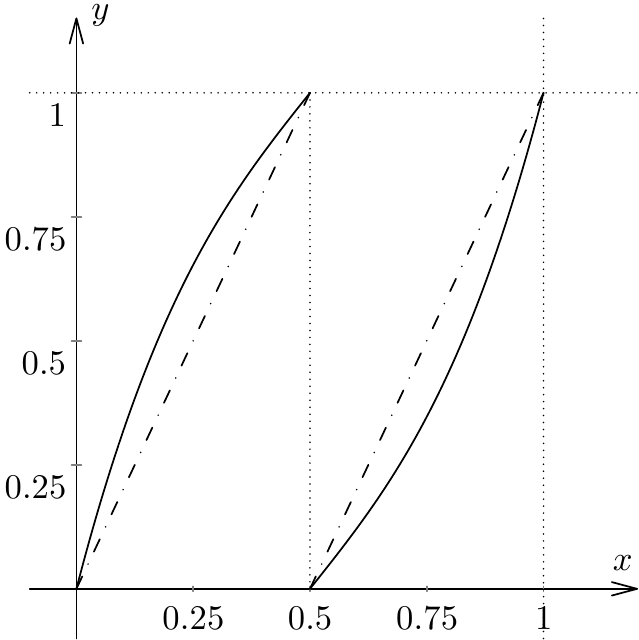}   }
 \caption{The solid curve is a  plot of the bent baker's map.  For comparison,
 the dashed line is the plot of the doubling map.}
 \end{figure}

\noindent The inverse branches $f_1, f_2: I \to I$ are defined by 
\begin{align*}
f_1(x) &= \frac{-2 + 2\sqrt{6}+ 2^{2/3}\left(9x + \sqrt{-38 + 18 \sqrt{6} + 
81x^2}\right)^{2/3}}{2^{11/6} \sqrt{-\frac{38}{3} + 6^{3/2} + 27 x^2}
\left(9x + \sqrt{-38 + 18\sqrt{6} + 81x^2} \right)^{1/3}} \\
f_2(x) &= \frac{-2 + 2\sqrt{6}+ 2^{2/3}\left(-9 + 9x + \sqrt{43 + 18\sqrt{6} - 
162 x+ 81x^2} \right)^{2/3}}{2^{11/6} \sqrt{-\frac{43}{3} + 6^{3/2} -54 x + 27 x^2}
\left(-9x + \sqrt{43 + 18\sqrt{6} - 162 x + 81x^2} \right)^{1/3}}
\end{align*}
and we can associate the transfer operators~$\mathcal L_t$ for $t \in \mathbb R$
according to~\eqref{eq:trop}. We next want to choose the following parameters for the computation. 
First, we compute $m=129$ Chebyshev nodes with accuracy of~$512$ digits. Then we
fix $\varepsilon = 10^{-75}$ and compute two polynomials~$p$ and $q$ using
Chebyshev collocation method. Afterwards, we use working precision of~$256$ for the routines
{\tt MinValue} and {\tt MaxValue}. The calculation gives  
\begin{align*}
    \alpha &= -\log \left( \mbox{\tt MaxValue} \frac {\mathcal L_{1+\varepsilon} p}{p} \right)= 
    0.6494631493\, 2069852907\, 6505\ldots \times 10^{-75}; \quad \mbox{ and  }\\
    \beta &= \log \left( \mbox{\tt MinValue} \frac {\mathcal L_{1-\varepsilon} q}{q}\right) = 
    0.6494631493\, 2069852907\, 7088 \ldots \times 10^{-75}.  
\end{align*}
We obtain the value of the Lyapunov exponent 
$$
\lambda(f,\mu) = 0.6494631493\,2069852907\, 6 \pm 10^{-21}.
$$ 
This is consistent with, and improves on,  Froyland's
estimate of $\lambda(f,\mu) =  0.64946$. 
We see that in this case the accuracy is  less than in 
other examples we have considered so far. One cause is the character of the inverse
branches $f_1$ and $f_2$: the formulae implies that providing we know the value
of $x\in(0,1)$ with an error of $10^{-k}$, we have the value of $f_{1}(x)$ and $f_{2}(x)$ with
an error of $10^{-k/6}$. 

Another source of complication is the diminished
hyperbolicity. A straightforward calculation gives that
$f'(x) \geq  f'(\frac{1}{2}) = 2 + 2 \sqrt{\frac{2}{3}} - \sqrt{6} =1.1835\ldots$.
This relatively weak hyperbolicity  also suggests  an explanation for why  the estimates are not as good as in the previous examples.   In particular, the maximal eigenfunction for $\mathcal L_t$ may be less regular (e.g.,  analytic on a relatively small Bernstein ellipse) which make the polynomial  approximation used in \S 4.2 less effective.

\section{Proof of  Theorem~\ref{thm:main}}
\label{s:proof}
In order to explain the proof of Theorem~\ref{thm:main} it helps to introduce
the following famous function from thermodynamic formalism. 
\subsection{Pressure function}
\label{s:pressure}
We begin by introducing the following well known definition.
\begin{definition}
To any map $f \in NC(I)$ we can associate the \emph{pressure function} 
$P\colon \mathbb R \to \mathbb R$ defined by
$$
P(t) = \lim_{n\to +\infty } \frac{1}{n} \log \sum_{f^nx=x} |(f^n)'(x)|^{-t} \hbox{ for } t \in \mathbb R.
$$
\end{definition}
This is one of many equivalent definitions of the pressure~\cite{walters}. 
The usefulness of the pressure function to study the Lyapunov exponent is shown
by the following simple lemma, the first three parts of which are well-known.  
\begin{lemma}\label{pressure}
 The pressure function  has the following properties:
\begin{enumerate}
\item $P(1) = 0$;
\item $P$ is an analytic convex  function; 
\item We can write $\lambda(f, \mu) = - \frac{dP(t)}{dt}|_{t=1}$ and;
\item For any $\varepsilon > 0$ we can write 
$$
-\frac{P(1+\varepsilon)}{\varepsilon} \leq - \frac{dP(t)}{dt}\Bigl|_{t=1} \leq \frac{P(1-\varepsilon)}{\varepsilon}. 
$$
\end{enumerate}
\end{lemma}

\begin{figure}[h]
\centerline{\includegraphics{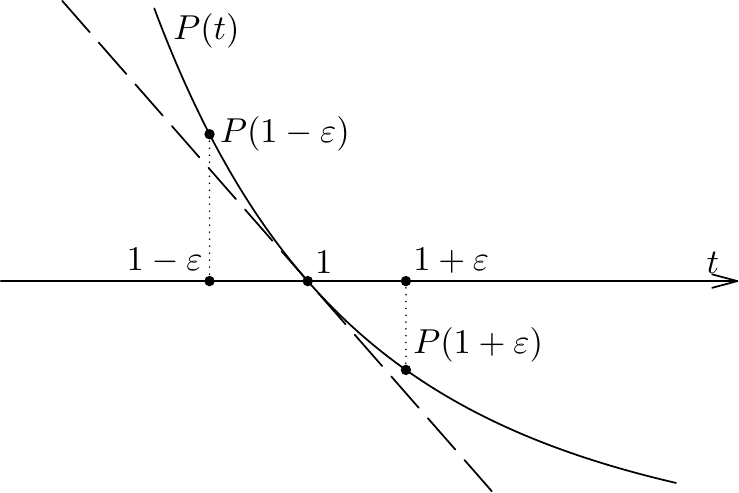}}
\caption{The pressure function $P(t)$ and the inequalities in part~(4) of the Lemma.}
\label{fig:pressure}
\end{figure}

\begin{proof}
The first three  parts are essentially due to Ruelle~\cite{ruelle} (see Corollary 5.27 and Exercise 5 (a) on p.99).  
The last observation follows easily from the convexity (see
Figure~\ref{fig:pressure}).
\end{proof}

This leads to the following useful bound on the Lyapunov exponent.  
\begin{cor}
    \label{cor:bounds}
For any $\varepsilon > 0$ the following double inequality holds 
$$  
-\frac{P(1+\varepsilon)}{\varepsilon}  \leq  \lambda(f, \mu) \leq \frac{P(1-\varepsilon)}{\varepsilon}.
$$
\end{cor}

\begin{proof}
This comes by substituting the identity in part 3 of Lemma \ref{pressure} into the inequality in part 4.
\end{proof}

\begin{rem}
At first sight, it may not seem very promising as an approach to estimating
$\lambda(f,\mu)$ to have to compute the pressures $P(1\pm \varepsilon)$ with an error $O(\varepsilon^2)$ in order 
to have an estimate on $\lambda(f,\mu)$ with error $O(\varepsilon)$. This means
that one has to estimate $P(1\pm \varepsilon)$ with the double accuracy of the
desired estimate for the Lyapunov exponent. 
Nevertheless it turns out that this approach is quite practical since it is quite easy to
estimate the pressure to high precision.  
\end{rem}

\subsection{Transfer operator for interval maps}
\label{ss:TrOp}
For definiteness, let us choose coordinates such that $I = [-1,1]$
(i.e., $a=-1$ and $b=1$ in Definition \ref{def:markovMap})
after a simple change of coordinates.  In addition, we shall also assume for simplicity that the map $f$ is
full-branch, i.e., $f((x_k, x_{k+1})) = (-1,1)$,
with inverse branches $f_k: (-1,1) \to (x_k, x_{k+1})$,
for $k=1, \cdots, n$, the general case being similar.


The approach to estimating the pressure is based on its interpretation in terms of the  family of 
transfer operators introduced in the Introduction.  These operators act  on the
Banach space~$\mathcal B$ of bounded analytic functions on domain $U_\rho \supset I$
enclosed by the Bernstein ellipse with the foci at~$0$ and~$1$ and given by
 $$
 \partial U_\rho = \left\{ z = \frac{1}{2} \left(\rho e^{i \theta} + \frac{e^{-i\theta}}{\rho}\right)  
  \colon  0 \leq \theta < 2\pi  \right\},\quad \rho > 1, 
 $$
 We define the norm on~$\mathcal B$ by $\|f\| = \sup_{z\in U} |f(z)| $. 
 In particular, by choosing~$\rho$ sufficiently close to~$1$ we can assume that
 the inverse branches~$f_j$ ($j=1, \cdots, n$) of the map $f \in NC(I)$ have
 analytic extensions to~$U_\rho$. Therefore the maps  $f_j: U_\rho \to U_\rho$ are well defined and their derivatives
 are non-zero, furthermore, since all~$f_j$ are contractions, we have that $\overline{\cup_j f_j U_\rho} \subsetneq U_\rho$. 
 We formally extend the definition of   the transfer operators from (2) as follows: 

\begin{definition}\label{transferoperator}
    A family of transfer operators $\mathcal L_t: \mathcal B  \to \mathcal B$ associated to
    $f \in NC(I)$ is defined by
      \begin{equation}
        \label{eq:tropa}
[\mathcal L_t h](x) = \sum_{k=1}^n |f_{j}^\prime(x)|^t h(f_{j}(x))
\qquad \mbox{ for } x\in U_\rho, \mbox{ and } t \in \mathbb R. 
\end{equation}
\end{definition}

\begin{rem}
    In this definition the functions  $|f_j'|^{t}$ are real valued and real
    analytic on $[-1,1]$. Thus by a slight abuse of notation we
    interpret~$|f_j'|^{t} $ as being the complex analytic extension of these functions 
    to~$U_\rho$.
\end{rem}

We can estimate the pressure values $P(t)$  using the maximal eigenvalue 
for the transfer operator $\mathcal L_t$.
\begin{lemma}
    \label{lem:spectral} 
    Let $f \in NC(I)$ and let~$\mathcal L_t \colon \mathcal B
    \to \mathcal B$ be the transfer operator defined by~\eqref{eq:tropa}. Then 
\begin{enumerate}
\item The spectral radius of $\mathcal L_t$ is $e^{P(t)}$.  
\item The rest of the spectrum is contained in a disk of radius strictly smaller
    than~$e^{P(t)}$.  
\item For any $h \in \mathcal B$ for which the restriction to $I$ is strictly positive   and any $x \in I$ we have 
    $\lim_{n \to+\infty} \left(\mathcal L_t^n h(x) \right)^{\frac{1}{n}} = e^{P(t)}$.
\end{enumerate}
\end{lemma}
\begin{proof}
Parts 1 and 2 are essentially due to Ruelle~\cite{ruelle} (see Proposition 5.13 and 5.24).
Part 3 follows directly from Part 1 and Part 2 and the classical spectral radius theorem (cf. ~\cite{ruelle}, Proposition 5.13 and 5.14).
\end{proof}
We can use this lemma to estimate the pressure values~$P(1\pm \varepsilon)$  
in Corollary~\ref{cor:bounds}. In particular, in order to estimate~$e^{P(t)}$
for $t \in \mathbb R$ we will use the following simple result.
\begin{lemma}
    \label{lem:tech}
    Assume that for~$t \in \mathbb R$     there exist a  function~$p \in
    \mathcal B$, strictly positive on~$I$, and a constant $\rho \in \mathbb R$ such that 
    $\sup_{x \in I}\frac{ \mathcal L_{t} p(x)}{p(x) }  \leq  e^{\rho} $
    then $ e^{P(t)} \le e^{\rho}$. 
\end{lemma}
\begin{proof}
Since $\mathcal L_{t} p(x) \leq  e^{\rho} p(x)$ for all $x \in I$ we can deduce that 
$\mathcal L_{t}^n p(x) \leq e^{n \rho} p(x)$ for $n \geq 1$.  Thus by Part 3 of
Lemma~\ref{lem:spectral} we have $e^{P(t)} = \lim_{n\to+\infty}
\left(\mathcal L_{t}^n p (x)  \right)^\frac{1}{n} 
 \leq e^\rho$ for any $x\in I$.
 \end{proof}
We now combine Lemma~\ref{lem:tech} and Corollary~3.3 to prove 
Theorem~\ref{thm:main}.
\begin{proof}[Proof of Theorem~\ref{thm:main}]
    By assumption, we know that there exist $0 < \alpha  <  \beta$ and two
    positive functions $p$ and $q$ such that 
    $$
    \sup_I \frac{\mathcal L_{1+\varepsilon} p}p \leq  e^{-\alpha}  \hbox{ and } \sup_I
    \frac{\mathcal L_{1-\varepsilon} q}q \leq  e^{\beta}. 
    $$
Applying   Lemma~\ref{lem:tech}  with  $t=1+\varepsilon$ and $\rho=-\alpha$
   then   gives
    $e^{P(1+\varepsilon)} \leq  e^{-\alpha}$ and thus $P(1+\varepsilon) \leq   -\alpha$, or equivalently, $-P(1+\varepsilon) \geq  \alpha$.
      On the other hand, 
    applying   Lemma~\ref{lem:tech}  with  $t=1-\varepsilon$  and $\rho=\beta$
   then   gives
    $e^{P(1-\varepsilon)} \leq e^\beta$ and thus  $P(1-\varepsilon) \leq  \beta$.
    Combining these two inequalities with Corollary~\ref{cor:bounds} we get 
    $- \frac{ \alpha}\varepsilon\leq  \lambda(f, \mu) \leq  \frac{
    \beta}\varepsilon $, as required.
\end{proof}

\begin{rem}
The above arguments extend easily to all maps $f \in NC(I)$, not necessary full
branch.  In particular, instead of a single domain~$U_\rho$, we consider  the disjoint 
union $U = \coprod_{k=1}^n U_\rho^{(k)}$
of domains 
 $U_\rho^{(k)} \supset [x_k, x_{k+1}]$ each bounded by a Bernstein ellipse with foci $x_n$ and $x_{n+1}$ 
 ($k=1, \cdots, n$).   The Banach space is now taken to be $\mathcal B = \oplus_{k=1}^n \mathcal B^{(k)}$
where $\mathcal B^{(k)}$ is the space of bounded analytic functions on $U_\rho^{(k)}$.  Finally, we use the extension of (1) 
to define $\mathcal L_t: \mathcal B \to \mathcal B$ by
$$
 [\mathcal L_t \underline h]_j(x) = \sum_{k: f((x_k,x_{k+1})) \cap (x_j,x_{j+1}) \neq
\varnothing} |f_{jk}^\prime(x)|^t h_k(f_{jk}(x))
 \qquad t \in \mathbb R, x\in U_\rho^{(j)}. 
$$
where $\underline h = (h_1, \cdots, h_n) \in \mathcal B$.  The argument then
proceeds as above.

\end{rem}
\section{Practical realisation}
\label{s:practice}
We next want to explain how to apply Theorem~\ref{thm:main} in practice. Below
we give one way of constructing test functions~$p$ and~$q$ that we used in order
to obtain estimates in the examples we considered. It is based on the
spectral Chebyshev collocation method. There are other methods  one
might consider, such as spline interpolation methods, proposed by Falk and
Nussbaum~\cite{FN}, but this approach suffices for our needs.

\subsection{Constructing test functions $p$ and $q$}
\label{ss:pq}

For notational simplicity, we will describe our construction in the special case of full branch maps.
The generalization to the general case of Markov maps
is fairly  straightforward where $I$ is replaced by the disjoint union of intervals. 

\begin{definition}
        \label{def:lagrpoly}
    Let $x_1 < x_2 < \ldots <  x_n$ be a collection of $n$ distinct real numbers. The Lagrange
    polynomials associated to $\{x_j\}_{j=1}^n$ are the polynomials
    \begin{equation}
        \label{eq:lagrpoly}
        \ell p_j(x): = \prod_{k\ne j} \frac{x-x_k}{x_j-x_k}, \qquad j = 1, \ldots, m. 
    \end{equation}
\end{definition}
 The Lagrange
    polynomials have the property that $\ell p_j(x_k) = \delta_j^k$ for all $j = 1, \ldots,
m$ and $k = 1, \ldots, n$. 
In a special case when the points $\{x_j\}_{j=1}^n$ are the roots of a certain
polynomial $p_n$, they can be written as 
\begin{equation}
    \label{eq:lagrpoly2}
    \ell p_j(x) = p_n(x) \cdot (p_n^\prime(x_j))^{-1} \cdot (x-x_j)^{-1}, \qquad
    j = 1, \ldots, m.
\end{equation}

We assume below that $I = [-1,1]$, the general case being similar after a simple
change of coordinates,  and $f \in NC(I)$. Let us assume that one wishes to
compute the Lyapunov exponent with an error of~$\delta$, in other words, we assume that one 
wishes to find an interval $(\lambda_1,\lambda_2) \ni \lambda(f,\mu)$ such that
$|\lambda_2 - \lambda_1| \le \delta$. In order to define the functions~$p$
and~$q$, we begin by choosing a natural number $m = m(\delta)$. Then we calculate numerically, with help
 of a computer, the following objects: 
\begin{enumerate}
\item Chebyshev nodes~$x_k: = \cos\left(\frac{\pi(2k+1)}{2m}\right) \in (-1,1)$, for $k =
    0,\ldots,m-1$ --- these are the roots of the Chebyshev polynomial of the
    first kind~$T_m$.
    In a general case of $I=[a,b]$ the Chebyshev
nodes have to be rescaled and shifted to~$I$ using the transformation $x\mapsto
x(b-a)+a$. 
    
The cosine function can be evaluated at a given point with
    arbitrary precision. In particular, in each of the Examples we consider
    we specify the number of digits~$N = N(\delta)$ requested in the actual program code.  
\item For $t = 1 \pm \varepsilon$ the matrices~$M^t \in GL(m,\mathbb R)$ given by 
    \begin{equation}
        \label{eq:matMt}
    M^t_{jk} :=[\mathcal L_t \ell p_j](x_k) =  \sum_{i = 1}^m |f_i^\prime(x_k)
    |^{t} \cdot  (\ell   p_j( f_i(x_k));
    \end{equation}
    Here we use the formula~\eqref{eq:lagrpoly2} to evaluate $\ell   p_j(
    f_i(x_k))$, using an inbuilt routine for evaluation of Chebyshev polynomials,
    which has guaranteed accuracy. 
\item The leading left eigenvectors~$v^t$ corresponding to the maximal eigenvalue
    of the matrices~$M^t$ for $t = 1 \pm \varepsilon$. They can be efficiently
    computed using the power method. 
\item 
    The polynomials $p$ and $q$ then given by linear combinations of Lagrange
polynomials $\ell p_j$ with coefficients coming from the eigenvectors:
\begin{equation}
    \label{eq:pq}
    p(x) = \sum_{j=0}^{m-1} v_j^{(1-\varepsilon)} \ell p_j(x) 
    \qquad \mbox{ and } \qquad 
    q(x) = \sum_{j=0}^{m-1} v_j^{(1+\varepsilon)} \ell p_j(x).
\end{equation}
However, the formula~\eqref{eq:pq} is prone to numerical errors. 
The polynomials~$p$ and~$q$ can also be written as a linear combination of Chebyshev
polynomials, and this has the advantage of being 
more computationally stable than the more direct expansion in terms of Lagrange
polynomials above. 
More precisely, the following expansion is well known. 
$$
p(x) = \sum_{j=0}^{m-1} a_j T_j(x), \qquad \mbox{ and }
q(x) = \sum_{j=0}^{m-1} b_j T_j(x), 
$$
where the coefficients are given in terms of the eigenvectors
$v^{(1-\varepsilon)}$ and $v^{(1+\varepsilon)}$: 
\begin{align*}
    a_j &= \frac2m \sum_{k=1}^m v_k^{(1-\varepsilon)} T_j(x_k) 
    &&\mbox{ and } \quad 
     b_j  = \frac2m \sum_{k=1}^m v_k^{(1+\varepsilon)} T_j(x_k), \mbox{ for
    } j = 1, \dots, m-1; \\ 
    a_0 & = \frac1m \sum_{k=1}^m v_k^{(1-\varepsilon)} 
    && \mbox{ and } \quad 
     b_0  = \frac1m \sum_{k=1}^m v_k^{(1+\varepsilon)}.
\end{align*}
This allows us to evaluate~$p$ and~$q$ efficiently. 
\item The supremums of the ratios $\frac{\mathcal L_{1+\varepsilon} p}p$ and
    $\frac{\mathcal L_{1-\varepsilon} q}q$ over the interval~$I$ is computed
    using internal routine {\tt MaxValue} with working precision set
    to~$D=D(\delta)$  digits. 
\end{enumerate}

In addition to exploiting the internal routine {\tt MaxValue} we can also 
apply Monte Carlo type method in order to carry out a heuristic check on its output. 

In the next  subsection we show that in the setting of uniformly expanding piecewise
analytic Markov maps the polynomials~$p$ and~$q$ satisfying the hypothesis of
Theorem~\ref{thm:main} can always be constructed, and moreover the conclusion of Theorem  \ref{thm:justification} holds.  

\subsection[Justification of the method]{Justification of the method: Proof of Theorem~\ref{thm:justification}}
\label{ss:justify}

We can denote by $\mathcal P_m \subset \mathcal B$   the polynomials on 
$I = [-1,1]$, say,  of degree~$m$.
We let $\pi_m: \mathcal B  \to \mathcal P_m$ be the projection onto the
polynomials of degree~$m$ given by the 
Chebychev--Lagrange collocation formula
$$
\pi_m(f)(x) = \sum_{j=0}^m f(x_j) \ell p_j(x), x \in I,
$$
where $x_j$, $j = 0, \ldots, m$ are the roots of the Chebyshev
polynomial~$T_{m+1}$, i.e. the Chebyshev nodes and~$\ell p_j$ ($j=0, \cdots, m$) are the
Lagrange polynomials on~$I$ associated to~$x_j$, $j = 0, \ldots, m$
defined by~\eqref{eq:lagrpoly}. 
In particular, we see that the restriction $\pi_m|_{\mathcal P_m}$ of~$\pi_m$
to~$\mathcal P_m$ is the identity.

The transfer operator
$\mathcal  L_{t}:\mathcal B \to \mathcal B$ defined by~\eqref{eq:tropa}
 is  compact, even nuclear, although this will not be needed.  
We require an  estimate on    the operator norm of the difference $ \mathcal L_t
- \mathcal L_t \pi_m$  defined  by
$$
\| \mathcal L_t - \mathcal L_t \pi_m\| = \sup_{\|f\|_{\mathcal B}=1}
\| (\mathcal L_t - \mathcal L_t \pi_m) (f) \|_{\mathcal B}.
$$

\begin{lemma}[see \cite{bs}, Theorem 3.3]\label{approx}
Let $f \in NC(I)$ and let $\mathcal L_t$ be the associated transfer operator
defined by~\eqref{eq:tropa}. Then there exists $C  > 0$ and $0 < \theta < 1$ such 
we can  bound that 
  $\| \mathcal L_t - \mathcal L_t \pi_m\|  \leq C \|\mathcal L_t\|\theta^m$ for $m \geq 1$.
\end{lemma}

This is also implicit in (\cite{wormell}, \S2.2).

\begin{rem}
    Although we cannot expect that~$\| I - \pi_m\|_{\mathcal B \to \mathcal B}
    \to 0 $ as $m\to \infty$, the composition with the  operator $\mathcal L_t$ 
allows the bound in Lemma~\ref{approx} since for any function~$f$ analytic
on~$U_\rho$ for some $\rho >1$, there exists $\rho^\prime > \rho$ such that   
the image~$\mathcal L_tf$ is analytic on~$U_{\rho^\prime}$. 
\end{rem}

It follows from the properties of~$\mathcal L_t$, Lemma~\ref{approx} and classical analytic 
perturbation (see the book of Kato~\cite{kato}) that we have the following:

\begin{lemma}\label{apert}
Let $f \in NC(I)$ and let $\mathcal L_t$ be the associated transfer operator
defined by~\eqref{eq:tropa}. Then for $\delta>0$ sufficiently small and~$m$ sufficiently large:
\begin{enumerate}
\item
 $\mathcal L_t \pi_m:  \mathcal B  \to \mathcal B$ has a
  simple maximal eigenvalue $\lambda_m$ with $|\lambda_m - e^{P(t)}| < \delta$;
 \item 
  The rest of the spectrum is contained in $\{z\in \mathbb C 
\colon |z|\ < e^{P(t)} - 2\delta\}$; and 
\item  
The corresponding eigenfunction~$h_m$ for~$\mathcal L_t \pi_m$ has a
restriction to~$I$ which  is strictly positive (i.e, $h_m(x) > 0$ for~$x\in I$).
\end{enumerate}
\end{lemma}
By  perturbation theory the positivity of the restriction of the
eigenfunction~$h$ associated to $e^{P(t)}$ for $\mathcal L_t$ onto~$I$ 
implies the same for~$h_m$ (since $\inf_{x\in I} |h(x) - h_m(x)| \leq \|h - h_m\|$ 
will be arbitrary small for $m$ sufficiently large).
The restriction $\pi_m \mathcal L_t|_{\mathcal P_m} $  is a finite rank operator 
 $\pi_m \mathcal L_t: \mathcal P_m \to \mathcal P_m$
given by 
  $$
\pi_m   \mathcal L_t :  g \mapsto \sum_{j=0}^m [\mathcal L_t g](x_j) \ell p_{j}, 
 $$
 where $x_j$, $j = 0, \ldots, m$ are the Chebyshev nodes introduced in
 Section~\ref{ss:pq}. Observe that in the basis of Lagrange polynomials $\{\ell
 p_j\}_{j=0}^m$ given by~\eqref{eq:lagrpoly}
 the operator~$\pi_m \mathcal L_t$ is given by the matrix~$M^t$ defined
 by~\eqref{eq:matMt}. 
In particular, maximal eigenvalue~$\lambda_m$ for~$\mathcal L_t\pi_m$ is also an
eigenvalue for the matrix~$M^t$ corresponding to the eigenvector 
$\pi_m(h_m) \in \mathcal P_m$.  This completes the proof of
Theorem~\ref{thm:justification}. 

\phantom{ping} \hfill $\square$

Given $f \in NC(I)$ and $\delta>0$ there exist $N = N(\delta)$ and  $D =
D(\delta)$ such that the polynomials~$p$ and~$q$ of degree $m \geq N$ 
with coefficients computed to~$D$ decimal places lead to estimates with
$\frac1\varepsilon\log\frac\beta\alpha < \delta$.  In particular, the
exponential convergence in Lemma~\ref{approx}
implies that $N(\delta) = O(\log(\delta/\varepsilon))$ 
and~$D(\delta) = O(-\log \delta, -\log \varepsilon)$. 
    
\begin{rem}[Heuristic estimates on the accuracy of approximation]
For small~$|\varepsilon|\ll1$ an~$O(\varepsilon^2)$ approximation to the
eigenvalue~$e^{P(t)}$ should give an~$O(\varepsilon)$ estimate on the Lyapunov exponent
(since we divide out by $\varepsilon$ in the formulae).  
   Furthermore,  this  error  is related to the (uniform) approximation error of
   the associated eigenfunction~$f$  by the interpolating polynomial based
   on~$m$  points, say, which is well known to be bounded by~$\|f^{k}\|_\infty(k+1)!$.  
   Even for very regular (e.g., analytic) functions~$f$ one only 
   expects~$\|f^{k}\|_\infty$ to tend to zero at best   exponentially fast. 
  Therefore, we might want~$\varepsilon \sim \sqrt{1/(k+1)!}$.
  In particular, for degree $k=10$ one gets $\varepsilon = 5.10^{-4}$,  
  for $k=20$ one gets $\varepsilon=6.10^{-10}$, and 
  for $k=100$ one gets $\varepsilon=1.10^{-79}$. 
\end{rem}


\medskip
      
  {\footnotesize
  \noindent
  \textsc{P. Vytnova, Department of Mathematics, Warwick University, Coventry,
  CV4 7AL, UK}
  \noindent
  \textit{E-mail address}: \texttt{P.Vytnova@warwick.ac.uk}
  \par
  \addvspace{\medskipamount}
  \noindent
  \textsc{M. Pollicott, Department of Mathematics, Warwick University, Coventry, CV4 7AL, UK.}
  \noindent
  \textit{E-mail address}: \texttt{masdbl@warwick.ac.uk}
}

\end{document}